\documentclass{amsart}
\usepackage{setspace, amsmath, amsthm, amssymb, amsfonts, amscd, epic, graphicx, ulem, dsfont}
\usepackage[T1]{fontenc}
\usepackage{multirow}
\usepackage{bbm}
\usepackage{enumerate}

\makeatletter \@namedef{subjclassname@2010}{
  \textup{2020} Mathematics Subject Classification}
\makeatother

\newtheorem{thm}{Theorem}[section]

\theoremstyle{remark}

\theoremstyle{definition}

\newcommand{\R}{\mathbb{R}}

\begin{document}

\title[Paranormal not normal]{A closed densely defined operator $T$ such that both $T$ and $T^*$ are injective and paranormal yet $T$ is not normal}
\author[M. H. MORTAD]{Mohammed Hichem Mortad}

\thanks{}
\date{}


\address{Department of
Mathematics, University of Oran 1, Ahmed Ben Bella, B.P. 1524, El
Menouar, Oran 31000, Algeria.\newline {\bf Mailing address}:
\newline Pr Mohammed Hichem Mortad \newline BP 7085 Seddikia Oran
\newline 31013 \newline Algeria}

\email{mhmortad@gmail.com, mortad.hichem@univ-oran1.dz.}

\begin{abstract}In this note, we give an example of a densely defined
closed one-to-one paranormal operator $T$ whose adjoint is also
injective and paranormal, but $T$ fails to be normal.
\end{abstract}

\maketitle

\section{Introduction}

To read this note easily, readers must have knowledge of linear
bounded and unbounded operators, as well as matrices of unbounded
operators. Some useful references related to the topics of the paper
are \cite{Kubrusly-book-operatopr-exercise-sol},
\cite{SCHMUDG-book-2012} and \cite{tretetr-book-BLOCK} respectively.

We recall a few definitions though: A linear operator $A$ with a
domain $D(A)$ contained in a Hilbert space $H$ is said to be densely
defined if $\overline{D(A)}=H$. Say that a linear operator $A$ is
closed if its graph is closed in $H\oplus H$. A densely defined
linear operator $A$ is called normal if $A$ is closed and
$AA^*=A^*A$. Equivalently, $\|Ax\|=\|A^*x\|$ for all $x\in
D(A)=D(A^*)$.

Recall also that a linear operator $A:D(A)\subset H\to H$ (where $H$
is a Hilbert space) is said to be \textbf{paranormal} if
\[\|Ax\|^2\leq \left\|A^2x\right\|\|x\|\]
for all $x\in D(A^2)$. This is clearly equivalent to $\|Ax\|^2\leq
\|A^2x\|$ for all \textit{unit} vectors $x\in D(A^2)$.

T. Ando \cite{Ando-Paranormal-tensor-product-sums et al} showed that
if $T\in B(H)$ is such that both $T$ and $T^*$ are paranormal and if
$\ker T=\ker T^*$, then $T$ is normal. The proof is not that obvious
as it is for the stronger class of hyponormal operators (where the
assumption $\ker T=\ker T^*$ is not even needed). 

In this paper, we show that the natural generalization to unbounded
closed operators is untrue even when both $T$ and $T^*$ are
injective. The main idea comes from an explicit example of a closed
densely defined $T$ such that (obtained recently in
\cite{Dehimi-Mortad-CHERNOFF}):
\[D(T^2)=D({T^*}^2)=\{0\}\]
(other examples of higher powers may be found in
\cite{Mortad-TRIVIALITY POWERS DOMAINS}). To digress, readers may
think a priori that the example here is weaker than Chernoff's
famous example \cite{CH}, given that the operator there is also
symmetric and semi-bounded. In fact, the two examples are just of
different calibers. Let us elaborate a little more. In Chernoff's
case, it is impossible to have $D({T^*}^2)=\{0\}$. Indeed, since $T$
is symmetric and densely defined, i.e. $T\subset T^*$, then
$T^*T\subset {T^*}^2$. By the closedness of $T$, it results that
${T^*}^2$ must be densely defined as $T^*T$ is self-adjoint, and in
particular densely defined. Another nuance, both $T$ and $T^*$ in
the example in \cite{Dehimi-Mortad-CHERNOFF} are also injective
(which was missed by the authors there).

\section{Main Counterexample}

The counterexample may be considered simple, however, it was
obtained from two recent papers as well as the powerful tool of
matrices of operators.

\begin{thm}On some Hilbert space, there is a closed densely defined
operator $T$ such that both $T$ and $T^*$ are one-to-one and
paranormal, yet $T$ is not normal.
\end{thm}

\begin{proof}The Hilbert space in question is
 $L^2(\R)\oplus L^2(\R)$. From \cite{Dehimi-Mortad-CHERNOFF}, we have an explicit example of a densely defined unbounded closed operator $T$ for
which:
\[D(T^2)=D({T^*}^2)=\{0\}.\]
More precisely,
\[T=\left(
      \begin{array}{cc}
        0 & A^{-1} \\
        B & 0 \\
      \end{array}
    \right)
\]
on $D(T):=D(B)\oplus D(A^{-1})\subset L^2(\R)\oplus L^2(\R)$, and
where $A$ and $B$ are two unbounded self-adjoint operators such that
\[D(A)\cap D(B)=D(A^{-1})\cap D(B^{-1})=\{0\}\]
where $A^{-1}$ and $B^{-1}$ are not bounded (as in \cite{KOS}).
Hence
\[T^*=\left(
      \begin{array}{cc}
        0 & B \\
        A^{-1}& 0 \\
      \end{array}
    \right)\]
    for $A^{-1}$ and $B$ are both self-adjoint. Observe now that
    both $T$ and $T^*$ are one-to-one since both $A^{-1}$ and $B$
    are so.

Both $T$ and $T^*$ are trivially paranormal because
$D(T^2)=D({T^*}^2)=\{0\}$. So paranormality of both operator need
only be checked at the zero vector and this is plain as
\[\|Tx\|^2=\|T^2x\|\|x\|=0\text{ and }\|T^*x\|^2=\|T^{*2}x\|\|x\|=0\]
for $x=0$.

However, $T$ cannot be normal for it were, $T^2$ would too be
normal, in particular it would be densely defined which is
impossible here.
\end{proof}

\section{An open problem}
Recall that thanks to the closed graph theorem, everywhere defined
closable operators are automatically bounded. This includes the
classes of symmetric and hyponormal operators among others. Since
some densely defined paranormal operators are not necessarily
closable as was already observed in
\cite{Daniluk-paranormals-non-closable} or in
\cite{Mortad-paranormal-daniluk}, it would be interesting to find an
unbounded paranormal operator which is everywhere defined. Can we
find such example?

\end{document}